%% file: arxiv.tex
\documentclass{article}

\usepackage[utf8]{inputenc}
\usepackage{authblk} % Authors
\usepackage{times}
\usepackage{amsthm,amsmath,amsfonts,amssymb,mathtools}
\usepackage[numbers]{natbib}
\usepackage{enumitem}
\usepackage[a4paper, textheight=21.9cm, textwidth=15.90cm]{geometry}
\usepackage{xcolor} % colored text
\usepackage[colorlinks,citecolor=blue, linkcolor = blue,urlcolor=blue]{hyperref}

\newif\ifarxiv\arxivtrue

\newcommand{\R}[0]{\mathbb{R}}
\newcommand{\Rd}[0]{\mathbb{R}^{d}}
\newcommand{\Rn}[0]{\mathbb{R}^{n}}
\newcommand{\Rplus}[0]{\mathbb{R}_{+}}
\newcommand{\natnum}[0]{\mathbb{N}^{*}}
\newcommand{\euler}{\operatorname{e}}
\newcommand{\iu}{{i\mkern1mu}}
\DeclareMathOperator*{\argmax}{arg\,max}
\newcommand{\E}[0]{\mathbb{E}}

\newtheorem{theorem}{Theorem}[section]

\newtheorem{lemma}[theorem]{Lemma}
\newtheorem{corollary}[theorem]{Corollary}
\newtheorem{remark}{Remark}[section]

\newtheorem{example}{Example}[section]

\newtheorem{assumption}{Assumption}[section]

\newcommand\extrafootertext[1]{%
    \bgroup
    \renewcommand\thefootnote{\fnsymbol{footnote}}%
    \renewcommand\thempfootnote{\fnsymbol{mpfootnote}}%
    \footnotetext[0]{#1}%
    \egroup
}

\title{On the orthogonality of zero-mean Gaussian measures: Sufficiently dense sampling}

\author{Reinhard Furrer}
\author{Michael Hediger}

\affil{Department of Mathematical Modeling and Machine Learning and Department of Mathematics, Winterthurerstrasse 190, 8057 Zurich, Switzerland}

\begin{document}

\maketitle

\extrafootertext{\textit{Email addresses:} \texttt{reinhard.furrer@uzh.ch} (R.\ Furrer), \texttt{michael.hediger@math.uzh.ch} (M.\ Hediger)}

\begin{abstract}
For a stationary random function $\xi$, sampled on a subset $D$ of $\mathbb{R}^{d}$, we examine the equivalence and orthogonality of two zero-mean Gaussian measures $\mathbb{P}_{1}$ and~$\mathbb{P}_{2}$ associated with $\xi$. We give the isotropic analog to the result that the equivalence of $\mathbb{P}_{1}$ and~$\mathbb{P}_{2}$ is linked with the existence of a square-integrable extension of the difference between the covariance functions of $\mathbb{P}_{1}$ and~$\mathbb{P}_{2}$ from $D$ to $\mathbb{R}^{d}$. We show that the orthogonality of $\mathbb{P}_{1}$ and~$\mathbb{P}_{2}$ can be recovered when the set of distances from points of $D$ to the origin is dense in the set of non-negative real numbers.\\[2mm]
\end{abstract}

\input{main.tex}

\section*{Acknowledgments}
The authors are grateful to Thomas Leh\'{e}ricy for several helpful
discussions. This work was supported by the Swiss National Science Foundation SNSF-175529.

\bibliographystyle{plain} 
\bibliography{references}

\end{document}

%% file: main.tex
%%%% Main text entry area:
\section{Introduction}

\subsection{Primary notation} \label{sec:not}
We use the notation $\natnum \coloneqq \mathbb{N}\setminus \{0\}$ and $\Rplus \coloneqq [0,\infty)$ for the set of non-zero natural numbers and non-negative real numbers, respectively. Given $x \in \Rd$, the Euclidean norm of $x$ is identified with $\lVert x \rVert \coloneqq \sqrt{\langle x, x \rangle}$, where $\langle x, y \rangle \coloneqq x^{\mathrm{t}}y$, $x,y \in \Rd$, is the dot product on $\Rd$. 
An open ball of radius $r \in \Rplus$, centered at $x \in \Rd$, is denoted with $B_{r}(x)$. Finally, the Borel $\sigma$-algebra on $\Rd$ is written as $\mathfrak{B}(\Rd)$. 

\subsection{General framework} \label{sec:head}

Let $\mathbb{X}$ denote the space of real valued functions defined on $\Rd$, i.e., $\mathbb{X} \coloneqq \{s \colon s(x) \in \R, x \in \Rd\}$. $\mathbb{X}$ shall be equipped with the $\sigma$-algebra $\mathcal{U}$, the smallest $\sigma$-algebra which contains the algebra of sets
\begin{equation} \label{CylinderAlgebra}
    \mathcal{A} \coloneqq \bigcup_{n=1}^{\infty}\bigcup_{\substack{x_{i} \in \Rd \\ 1 \leq i \leq n}}\sigma\big(\{C_{x_{1}, \dotsc, x_{n}}(B_{n}) \colon B_{n} \in \mathfrak{B}(\Rn)\}\big),
\end{equation}
with $C_{x_{1}, \dotsc, x_{n}}(B_{n}) \coloneqq \{s \in \mathbb{X} \colon (s(x_{1}), \dotsc, s(x_{n})) \in B_n\}$, $B_n \in \mathfrak{B}(\Rn)$, being the cylinder sets on $\mathbb{X}$ over the coordinates $x_{1}, \dotsc, x_{n}$. The $\sigma$-algebra $\mathcal{U}$ is referred to as the $\sigma$-algebra generated by the cylinder sets. On a measure space $(\Omega, \mathcal{F})$ we consider a random function $\xi \colon \Omega \to \mathbb{X}$ which is $\mathcal{F}/\mathcal{U}$ measurable. That is, we consider a random field $\{\xi_{x} \colon x \in \Rd\}$ which has, for fixed $\omega \in \Omega$, real valued sample paths $\xi(\omega) = [x \mapsto \xi_{x}(\omega)]$ defined on $\Rd$. Let $\mathbb{P}_{1}$ and $\mathbb{P}_{2}$ be two probability measures, defined on $\mathcal{F}$, under which $\xi$ is Gaussian, i.e., for $\ell = 1,2$, $(\xi_{x_{1}}, \dotsc, \xi_{x_{n}})$ is a Gauss vector under $\mathbb{P}_{\ell}$ for any $n \in \natnum$ and $x_{1}, \dotsc, x_{n} \in \Rd$. In particular, $\mathbb{P}_{1}$ and $\mathbb{P}_{2}$ are said to be Gaussian measures on $\sigma(\xi) \coloneqq \{\xi^{-1}(U) \colon U \in \mathcal{U}\}$, the $\sigma$-algebra generated by $\xi$. Given a real valued random variable $Y$, defined on $(\Omega, \mathcal{F})$, we write $\E_{\ell}[Y] \coloneqq \int_{\Omega}Y(\omega)\mathbb{P}_{\ell}(d\omega)$ for the mean of $Y$ under $\mathbb{P}_{\ell}$, $\ell = 1,2$. The mean and covariance functions of $\xi$ under $\mathbb{P}_{\ell}$ are denoted with $\mu_{\ell}(x) \coloneqq \E_{\ell}[\xi_{x}]$ and $c_{\ell}(x,y) \coloneqq \E_{\ell}[(\xi_{x}-\mu_{\ell}(x))(\xi_{y}-\mu_{\ell}(y))]$, $\ell = 1,2$. Given a subset $D \subset \Rd$, we define the sub $\sigma$-algebra $\mathcal{U}_{D} \subset \mathcal{U}$ upon taking the union in \eqref{CylinderAlgebra} over $D$ instead of $\Rd$. We then denote $\sigma_{D}(\xi) \coloneqq \{\xi^{-1}(U), U \in \mathcal{U}_{D}\}$ as the sub $\sigma$-algebra ($\sigma_{D}(\xi) \subset \sigma(\xi)$) generated by $\xi$ with sample paths $x \mapsto \xi_{x}(\omega)$, $\omega \in \Omega$, which are restricted to $D$. We call $D$ the sampling domain for $\xi$. 

\subsection{On the equivalence and orthogonality of Gaussian measures with different covariance functions} \label{sec:EquivOrtho}

Conditions for the equivalence or orthogonality of the Gaussian measures $\mathbb{P}_{1}$ and $\mathbb{P}_{2}$ on $\sigma_{D}(\xi)$ are discussed in Chapter~III of \cite{GaussianRandomProcesses}. Other common references are \cite{GikhmanSkorokhod} (Chapter~VII) or also \cite{Yadrenko} (Chapter~III). Let us recall some terminology. $\mathbb{P}_{2}$ is said to be absolutely continuous with respect to $\mathbb{P}_{1}$ on $\sigma_{D}(\xi)$ if $\mathbb{P}_{1}(A) = 0$ implies $\mathbb{P}_{2}(A) = 0$, $A \in \sigma_{D}(\xi)$. $\mathbb{P}_{1}$ and $\mathbb{P}_{2}$ are termed equivalent if they are mutually absolutely continuous. On the other hand, $\mathbb{P}_{1}$ and $\mathbb{P}_{2}$ are referred to as orthogonal on $\sigma_{D}(\xi)$ if there exist $A \in \sigma_{D}(\xi)$ for which $\mathbb{P}_{1}(A) = 0$ but $\mathbb{P}_{2}(A) = 1$. Orthogonal measures are denoted by $\mathbb{P}_{1} \perp \mathbb{P}_{2}$. Using Lebesgue's decomposition theorem, it can be shown that $\mathbb{P}_{1} \perp \mathbb{P}_{2}$ on $\sigma_{D}(\xi)$ if there exists $(A_{n}) \subset \sigma_{D}(\xi)$ such that $\mathbb{P}_{1}(A_{n}) \xrightarrow[]{} 0$ but $\mathbb{P}_{2}(A_{n}) \xrightarrow[]{} 1$ as $n \xrightarrow[]{} \infty$ (compare to p.\ 64 in \cite{GaussianRandomProcesses}). Further, it is well known that Gaussian measures are either equivalent or orthogonal (see for instance Theorem~1 in Chapter~III of \cite{GaussianRandomProcesses}).

Throughout this article we assume that $\xi$ is stationary (homogeneous) under $\mathbb{P}_{1}$ and $\mathbb{P}_{2}$. That is, for $\ell = 1,2$, $\xi$ has constant mean function $\mu_{\ell}$ and covariance function $c_{\ell}$ that depends only on the difference $x-y$, $x,y \in \Rd$. Furthermore, the following two items are assumed to be satisfied: \begin{enumerate}[label=(\roman*)]
    \item \label{a1.1} $\mu_{\ell}(x) = 0$, $\ell = 1,2 \, ;$ 
    \item \label{a1.2} for $\ell = 1,2$, there exists a finite measure $F_{\ell}$, uniquely defined on $\mathfrak{B}(\Rd)$, such that $$c_{\ell}(x,y) = \int_{\Rd}\!\euler^{\iu \langle \lambda,x-y  \rangle}F_{\ell}(d\lambda).$$
\end{enumerate}
Recall that if for any $x_{0} \in \Rd$, $\mathbb{E}_{\ell}[\lvert\xi_{x}-\xi_{x_{0}}\rvert^{2}] \xrightarrow[]{} 0$ as $\lVert x-x_{0} \rVert \xrightarrow[]{} 0$, $\xi$ is said to be continuous in mean-square (m.s.\ continuous) under $\mathbb{P}_{\ell}$. We remark that under the assumption of stationarity, \ref{a1.2} is satisfied when $\xi$ is m.s.\ continuous under both $\mathbb{P}_{1}$ and $\mathbb{P}_{2}$ (see Theorem~2 in Section~2 of Chapter~IV in \cite{GikhmanSkorokhod}). Further, since for $\ell = 1,2$, $\xi$ is stationary under $\mathbb{P}_{\ell}$, there exists $k_{\ell} \colon \Rd \to \mathbb{R}$ such that for any $x,y \in \Rd$, $c_{\ell}(x,y) = k_{\ell}(x-y)$. Then, see p.\ 208 in the latter reference, $\xi$ is m.s.\ continuous under $\mathbb{P}_{\ell}$ if and only if $k_{\ell}$ is continuous at zero.

Chapter~III.4.2 of \cite{GaussianRandomProcesses} gives an overview of the case where Gaussian measures differ only in terms of their covariance function. In particular, the equivalence or orthogonality of $\mathbb{P}_{1}$ and $\mathbb{P}_{2}$ on $\sigma_{D}(\xi)$ is linked to the difference between the covariance functions $c_{1}$ and $c_{2}$ restricted to $D \times D$, 
\begin{equation*} \label{differenceKernels}
    \delta(x,y) \coloneqq c_{1}(x,y) - c_{2}(x,y), \quad x,y \in D.
\end{equation*}

For further reference, the following assumption is needed. 

\begin{assumption}[Absolute continuity of spectral measures]\label{a2} 
For $\ell = 1,2$, the spectral measure $F_{\ell}$ is absolutely continuous with respect to the Lebesgue measure on $\mathfrak{B}(\Rd)$ with spectral density $f_{\ell}(\lambda) = F_{\ell}(d\lambda)/d\lambda$.
\end{assumption}

\begin{remark} \label{remark1}
    It follows from \ref{a1.2} that for $\ell = 1,2$, the spectral density $f_{\ell}$ must be non-negative a.e.\ with respect to the Lebesgue measure on $\mathbb{R}^{d}$. To see it, we recall that because of \ref{a1.2} there exists a stochastic orthogonal measure $\zeta_{\ell}$, defined on $\mathfrak{B}(\Rd)$, with structure function $F_{\ell}$, s.t.\ for any $x \in \mathbb{R}^{d}$ we have that $\mathbb{P}_{\ell}$ a.s.,
    \begin{gather*}
        \xi_{x} = \int_{\mathbb{R}^{d}}\!\euler^{\iu \langle \lambda,x\rangle}\zeta_{\ell}(d\lambda).
    \end{gather*}
    See for instance Theorem~1 in Section~5 of Chapter~IV in \cite{GikhmanSkorokhod}. Then, if we let $A = f_{\ell}^{-1}\big((-\infty,0)\big)$, by Assumption~\ref{a2}, $f_{\ell}$ is Borel measurable and hence $A \in \mathfrak{B}(\Rd)$. In particular,
    \begin{gather*}
        \int_{A}\!f_{\ell}(\lambda)d\lambda = F_{\ell}(A) = \mathbb{E}_{\ell}\big[\lvert \zeta_{\ell}(A) \rvert^{2}\big],
    \end{gather*}
    which shows that $A$ must have Lebesgue measure zero.
\end{remark}

The following result serves as an anchor point for our study. For a proof, we refer to Section~III.4.2 in \cite{GaussianRandomProcesses} (see Theorem~11).  

\begin{theorem} \label{thm1}
    Suppose that Assumption~\ref{a2} is satisfied where $f_{1}$ and $f_{2}$ are bounded on $\Rd$. Then, the Gaussian measures $\mathbb{P}_{1}$ and $\mathbb{P}_{2}$ are equivalent on $\sigma_{D}(\xi)$ if and only if the restriction $\delta$ satisfies the following properties:
    \begin{enumerate}[label=(\alph*)]
    \item \label{thm1.1} There exists extension $\prescript{\scalebox{0.5}{$\uparrow$}}{}{}\delta$ of $\delta$ to $\Rd \times \Rd$ which is square-integrable, i.e., $$\int_{\Rd}\int_{\Rd}\lvert\prescript{\scalebox{0.5}{$\uparrow$}}{}{}\delta(x,y)\rvert^{2}dxdy < \infty \, ;$$
    \item \label{thm1.2} The Fourier transform $\varphi$ of $\prescript{\scalebox{0.5}{$\uparrow$}}{}{}\delta$ satisfies
    \begin{gather*}
        \int_{\Rd}\int_{\Rd}\frac{\lvert\varphi(\lambda,\mu)\rvert^{2}}{f_{1}(\lambda)f_{2}(\mu)}d\lambda d\mu < \infty.
    \end{gather*}
\end{enumerate}
\end{theorem}
Note that the proof given in \cite{GaussianRandomProcesses} is based on random functions that have sample paths defined on $\R$ instead of $\Rd$. Nevertheless, the arguments proposed for the case $d = 1$ can be recycled to prove the case where $d > 1$. We also remark that the existence of the spectral density $f_{\ell}$ is guaranteed if $k_{\ell}(x)$ is absolutely integrable on $\Rd$ (see for instance p.\ 211 in \cite{GikhmanSkorokhod}). 
\medskip

Theorem~\ref{thm1} allows us to easily deduce the orthogonality of $\mathbb{P}_{1}$ and $\mathbb{P}_{2}$ when $D$ is chosen to be the entire $\Rd$. In particular, under the assumption of bounded spectral densities, Theorem~\ref{thm1} shows that if $f_{1}$ and $f_{2}$ differ on a set of positive Lebesgue measure, $\mathbb{P}_{1}$ and $\mathbb{P}_{2}$ must be orthogonal on $\sigma(\xi)$ (see for instance p.\ 95 in \cite{GaussianRandomProcesses}). Using the two-dimensional Hankel transform (\cite{Hankel}) we will give the analog of Theorem~\ref{thm1} when $\xi$ is isotropic under $\mathbb{P}_{1}$ and $\mathbb{P}_{2}$ (Theorem~\ref{thm2}). Recall that $\xi$ is isotropic under $\mathbb{P}_{\ell}$ if $k_{\ell}(x)$ is a function of $\lVert x \rVert$ only. In a further step, we aim to recover the orthogonality of $\mathbb{P}_{1}$ and $\mathbb{P}_{2}$ when $D$ is dense in $\Rd$ (respectively $D_{+} \coloneqq \{\lVert x \rVert \colon x\in D\}$ is dense in $\Rplus$). Specifically, we will prove that if the covariance functions $c_{1}$ and $c_{2}$ are uniformly continuous on $\Rd \times \Rd$, the density of $D$ in $\Rd$ is enough to obtain orthogonal measures $\mathbb{P}_{1}$ and $\mathbb{P}_{2}$ (Theorem~\ref{thm3}) --- this under the assumption that $f_{1}$ and $f_{2}$ are bounded and different on a set of positive Lebesgue measure. If $\xi$ is isotropic, we show that one arrives at the same conclusion if $D_{+}$ is dense in $\Rplus$ (Theorem~\ref{thm4}). The latter result allows us to easily deduce the orthogonality of $\mathbb{P}_{1}$ and $\mathbb{P}_{2}$ when $\xi$ is sampled along a continuous and unbounded path which starts at zero. As an example, we will deduce the almost sure orthogonality of $\mathbb{P}_{1}$ and $\mathbb{P}_{2}$ when $\xi$ is sampled along a $d$-dimensional Brownian motion which starts at the origin (Example~\ref{exampleBM}). In particular, Theorems~\ref{thm3} and \ref{thm4} complement the results given in \cite{GaussianRandomProcesses} by deducing the orthogonality of two Gaussian measures based on a countable and unbounded collection of points sampled in $\Rd$. In terms of an illustration, we will revisit the relationship between orthogonal families of Gaussian distributions and covariance parameter estimation (\cite{Bachoc}) and discuss conditions under which consistent estimators can be obtained (Theorem~\ref{consistency}).  

\section{Main results} \label{sec:main_results}

\subsection{Continuous extension} \label{sec:main_results1}
We use Theorem~\ref{thm1} as a starting point and note that with regard to the necessity of the imposed conditions, a slight modification is possible. In particular, the extension $\prescript{\scalebox{0.5}{$\uparrow$}}{}{}\delta$ of item \ref{thm1.1} in Theorem~\ref{thm1} can be chosen to be continuous. This follows from the fact that the measures $F_{1}$ and $F_{2}$ are finite. More explicitly, from Theorem~8 in Section~III.3 of \cite{GaussianRandomProcesses} we can see that the Gaussian measures $\mathbb{P}_{1}$ and $\mathbb{P}_{2}$ are equivalent on $\sigma_{D}(\xi)$ if and only if the restriction $\delta$ permits a representation
    \begin{gather*}
    \delta(x,y) = \int_{\Rd}\int_{\Rd}\euler^{-\iu(\langle\lambda,x\rangle - \langle\mu,y\rangle)}\Psi(\lambda,\mu)F_{1}(d\lambda)F_{2}(d\mu),
    \end{gather*}
with $\Psi$ satisfying $\int_{\Rd}\int_{\Rd}\lvert \Psi(\lambda,\mu) \rvert^{2}F_{1}(d\lambda)F_{2}(d\mu) < \infty$. Then, since $F_{1}$ and $F_{2}$ are finite, by Hölder's inequality,
\begin{gather*}
    \int_{\Rd}\int_{\Rd}\lvert \Psi(\lambda,\mu) \rvert F_{1}(d\lambda)F_{2}(d\mu) < \infty.
\end{gather*}
Thus, since $f_{1}$ and $f_{2}$ are non-negative a.e.\ (see Remark~\ref{remark1}), the function $\Psi(\lambda, \mu)f_{1}(\lambda)f_{2}(\mu)$ is absolutely integrable on $\Rd \times \Rd$. Therefore, its Fourier transform $\prescript{\scalebox{0.5}{$\uparrow$}}{}{}\delta$ is continuous and absolutely integrable on $\Rd \times \Rd$. But, since $f_{1}$ and $f_{2}$ are also bounded, we have that $\prescript{\scalebox{0.5}{$\uparrow$}}{}{}\delta$ is also square-integrable on $\Rd \times \Rd$. Further, on $D \times D$, $\prescript{\scalebox{0.5}{$\uparrow$}}{}{}\delta$ agrees with $\delta$. Hence, we have proven the following result:

\begin{theorem} \label{thm1_extension}
    Suppose that Assumption~\ref{a2} is satisfied where $f_{1}$ and $f_{2}$ are bounded on $\Rd$. Then, if the Gaussian measures $\mathbb{P}_{1}$ and $\mathbb{P}_{2}$ are equivalent on $\sigma_{D}(\xi)$, there exists a continuous extension $\prescript{\scalebox{0.5}{$\uparrow$}}{}{}\delta$ of $\delta$ to the entire $\Rd \times \Rd$ which is absolutely and square-integrable on $\Rd \times \Rd$.  
\end{theorem}

We will see later that the above relation between equivalent measures and continuous extensions of $\delta$, will allow us to easily read the orthogonality of $\mathbb{P}_{1}$ and $\mathbb{P}_{2}$ from the uniform continuity of $c_{1}$ and $c_{2}$ if $D$ is dense in $\Rd$. Before we arrive there, we establish analogous versions of Theorems~\ref{thm1} and~\ref{thm1_extension}, when $\xi$ is not only stationary but also isotropic.

\subsection{Isotropic random fields} \label{sec:main_results2}

Given $x \in \Rd$, we adopt a polar coordinate system $x = (r_{x}, \theta_{x})$, $r_{x} \in \Rplus$, $\theta_{x} \in \mathbb{S}^{d-1}$, where $\mathbb{S}^{d-1}$ is the unit sphere in $\Rd$. On $L^{2}(\mathbb{S}^{d-1})$ we consider a real orthonormal basis composed of spherical (surface) harmonics $S_{m}^{l}$, $l = 1, \dotsc, h(m,d)$ of degree $m \in \mathbb{N}$ (see Chapter XI, Section~11.3 in \cite{Erdelyi} or also Chapter IV, Section~2 of \cite{WeissStein}). We recall that
\begin{gather*}
    h(m,d) = \begin{cases}
        1, & m = 0,\\
        d, & m = 1, \\
        \binom{d+m-1}{m} - \binom{d+m-3}{m-2}, & m \geq 2.
    \end{cases}
\end{gather*}

\begin{assumption} [Isotropy]\label{a3} 
$\xi$ is isotropic under $\mathbb{P}_{\ell}$, $\ell = 1,2$.
\end{assumption} 

Under Assumption~\ref{a3}, since \ref{a1.2} of Section~\ref{sec:EquivOrtho} is satisfied, we have, for $\ell = 1,2$, and any $x,y \in \Rd$ (see (4.145) of \cite{Yaglom}),
\begin{equation} \label{radial1}
    c_{\ell}(x,y) = K_{d}^{2}\sum_{m=0}^{\infty}\sum_{l=1}^{h(m,d)}\!\int_{0}^{\infty}S_{m}^{l}(\theta_{x})\frac{J_{m+\frac{d-2}{2}}(\kappa r_{x})}{(\kappa r_{x})^{\frac{d-2}{2}}}S_{m}^{l}(\theta_{y})\frac{J_{m+\frac{d-2}{2}}(\kappa r_{y})}{(\kappa r_{y})^{\frac{d-2}{2}}} \Phi_{\ell}(d \kappa),
\end{equation}
where $K_{d}^{2} = 2^{d-1}\Gamma(d/2)\pi^{d/2}$, $\Gamma$ is the Gamma function, $J_{m + (d-2)/2}$ is the Bessel function of the first kind of order $m + (d-2)/2$ and $\Phi_{\ell}$ is a finite measure on $\Rplus$ defined upon $\Phi_{\ell}([a,b)) = F_{\ell}(B_{b}(0)\setminus B_{a}(0))$.

We consider the real vector space of sequences of functions $a(\kappa) \coloneqq (a_{m}^{l}(\kappa))$, $\kappa \in \Rplus$, $l = 1, \dotsc, h(m,d)$, $m \in \mathbb{N}$. On the latter vector space, we introduce the inner product
\begin{gather*}
    \langle a, b \rangle_{\Phi_{\ell}} \coloneqq \sum_{m=0}^{\infty}\sum_{l=1}^{h(m,d)}\int_{0}^{\infty} a_{m}^{l}(\kappa)b_{m}^{l}(\kappa)\Phi_{\ell}(d \kappa), \quad \ell = 1,2.
\end{gather*}
Further, we define $\mathcal{L}_{D}^{0}$ as the linear span over $\R$ of the set of sequences of functions 
\begin{gather*}
    \bigg\{a \colon a(\kappa) = \bigg(K_{d}S_{m}^{l}(\theta_{x})\frac{J_{m+\frac{d-2}{2}}(\kappa r_{x})}{(\kappa r_{x})^{\frac{d-2}{2}}}\bigg),\ x = (r_{x}, \theta_{x}) \in D\bigg\},
\end{gather*}
and introduce the correspondence
\begin{equation} \label{coresp}
    \eta(a) \coloneqq \sum_{k=1}^{N}\beta_{k}\xi_{x_{k}}, \quad a = \sum_{k=1}^{N}\beta_{k}a_{k} \in \mathcal{L}_{D}^{0}.
\end{equation}
Let $L_{D}(\mathbb{P}_{\ell})$ be the linear span of $\{\xi_{x} \colon x \in D\}$ over $\R$ under $\mathbb{P}_{\ell}$, $\ell = 1,2$. Given $\ell = 1,2$, we view $L_{D}(\mathbb{P}_{\ell})$ as a subspace of the inner product space $L^{2}(\Omega, \mathcal{F}, \mathbb{P}_{\ell})$. Then, we readily see that for any $a,b \in \mathcal{L}_{D}^{0}$,
\begin{equation} \label{coresp1}
    \langle a, b \rangle_{\Phi_{\ell}} = \E_{\ell}[\eta(a) \eta(b)], \quad \ell = 1,2.
\end{equation}
Thus, for $\ell = 1,2$, \eqref{coresp} provides an isometric correspondence between the inner product space $\mathcal{L}_{D}^{0}$ (equipped with $\langle \cdot, \cdot \rangle_{\Phi_{\ell}}$) and $L_{D}(\mathbb{P}_{\ell})$. Note that if we define $\mathcal{L}_{D}(\Phi_{\ell})$, $\ell = 1,2$, as the closure of $\mathcal{L}_{D}^{0}$ with respect to $\langle \cdot, \cdot \rangle_{\Phi_{\ell}}$, the isometric correspondence \eqref{coresp} can be extended to an isometric correspondence between $\mathcal{L}_{D}(\Phi_{\ell})$ and the closure of $L_{D}(\mathbb{P}_{\ell})$.

Suppose that there exists $a \in \mathcal{L}_{D}^{0}$ such that $\lVert a \rVert_{\Phi_{1}} \neq 0$ but $\lVert a \rVert_{\Phi_{2}} = 0$. Then, the Gaussian measures $\mathbb{P}_{1}$ and $\mathbb{P}_{2}$ are orthogonal on $\sigma_{D}(\xi)$. This is because $\xi$ is Gaussian with a zero-mean function. Hence the latter assumption allows us to conclude that
\begin{gather*}
\mathbb{P}_{1}\big(\eta(a) = 0\big) = 0 \quad \text{but} \quad \mathbb{P}_{2}\big(\eta(a) = 0\big) = 1.
\end{gather*}
We recall that the two norms $\lVert \cdot \rVert_{\Phi_{1}}$ and $\lVert \cdot \rVert_{\Phi_{2}}$ are termed equivalent on $\mathcal{L}_{D}^{0}$ if
\begin{gather} \label{equivNorms}
    \exists \ C_{1}, C_{2} > 0 \ \text{s.t.\ }\forall \ a \in \mathcal{L}_{D}^{0},\ 0 < C_{1} \lVert a \rVert_{\Phi_{2}} \leq \lVert a \rVert_{\Phi_{1}} \leq C_{2} \lVert a \rVert_{\Phi_{2}} < \infty.
\end{gather}
Then, it can be shown that the Gaussian measures $\mathbb{P}_{1}$ and $\mathbb{P}_{2}$ are orthogonal on $\sigma_{D}(\xi)$ if the condition \eqref{equivNorms} is violated (compare to Section III.1.3 in \cite{GaussianRandomProcesses}). Notice that if \eqref{equivNorms} is satisfied, then by construction of $\mathcal{L}_{D}(\Phi_{1})$ and $\mathcal{L}_{D}(\Phi_{2})$, \eqref{equivNorms} remains true with $\mathcal{L}_{D}^{0}$ replaced with either $\mathcal{L}_{D}(\Phi_{1})$ or $\mathcal{L}_{D}(\Phi_{2})$. The following lemma is of central importance.

\begin{lemma} \label{intermediateLemma}
    Suppose that Assumption~\ref{a3} is satisfied. Then, the Gaussian measures $\mathbb{P}_{1}$ and $\mathbb{P}_{2}$ are equivalent on $\sigma_{D}(\xi)$ if and only if for any pair $m,i \in \mathbb{N}$, $l = 1, \dotsc, h(m,d)$, $j = 1, \dotsc, h(i,d)$, the difference 
    \begin{gather} \label{delta_mi}
    \delta_{m,i}^{l,j}(r_{x}, r_{y}) \coloneqq \int_{\mathbb{S}^{d-1}}\int_{\mathbb{S}^{d-1}}\delta(x,y)S_{m}^{l}(\theta_{x})S_{i}^{j}(\theta_{y})d\theta_{x}d\theta_{y}, \quad x,y \in D, 
    \end{gather}
is representable as
\begin{gather*}
    \delta_{m,i}^{l,j}(r_{x}, r_{y}) = K_{d}^{2}\int_{0}^{\infty}\int_{0}^{\infty}\frac{J_{m+\frac{d-2}{2}}(r_{x}\kappa)}{(r_{x}\kappa)^{\frac{d-2}{2}}}\frac{J_{i+\frac{d-2}{2}}( r_{y}\iota)}{(r_{y}\iota)^{\frac{d-2}{2}}}\psi_{m,i}^{l,j}(\kappa, \iota)\Phi_{1}(d \kappa)\Phi_{2}(d \iota),
\end{gather*}
with
\begin{equation} \label{squarepsi}
\sum_{m=0}^{\infty}\sum_{l=1}^{h(m,d)}\!\int_{0}^{\infty}\sum_{i=0}^{\infty}\sum_{j=1}^{h(i,d)}\!\int_{0}^{\infty}\lvert \psi_{m,i}^{l,j}(\kappa, \iota)\rvert^{2}\Phi_{1}(d \kappa)\Phi_{2}(d \iota) < \infty.
\end{equation}

\end{lemma}

\begin{proof}
We define the functions 
\begin{gather*}
    \Rd \times \Rd \ni (\lambda, \mu) \mapsto \mathfrak{e}_{x,y}(\lambda, \mu) \coloneqq \euler^{\iu (\langle\lambda, x \rangle - \langle\mu, y \rangle)}, \quad x,y \in D.
\end{gather*}
Then, as on p.\ 82 of \cite{GaussianRandomProcesses}, we let $L^{0}_{D \times D}$ denote the linear space of functions of the form
\begin{gather*}
    u(\lambda, \mu) = \sum_{p,q}\beta_{pq}\mathfrak{e}_{x_{p},y_{q}}(\lambda, \mu),
\end{gather*}
where $x_{p}, y_{q} \in D$ and $\beta_{pq}$ are real coefficients. Further, the Hilbert space $L_{D \times D}(F_{1} \times F_{2})$ is defined as the closure of $L^{0}_{D \times D}$ with respect to the inner product
\begin{gather*}
    \langle u , v \rangle_{F_{1} \times F_{2}} \coloneqq \int_{\Rd} \int_{\Rd} u(\lambda, \mu) \overline{v(\lambda, \mu)} F_{1}(d\lambda)F_{2}(d\mu).
\end{gather*}
Using polar coordinates, we introduce the class of sequences of functions
\begin{gather*}
    \Rplus \times \Rplus \ni (\kappa, \iota) \mapsto \mathfrak{a}_{x,y}(\kappa, \iota) \coloneqq \big(\mathfrak{a}_{m,i}^{l,j}(\kappa, \iota)\big), \quad x = (r_{x}, \theta_{x}), \ y = (r_{y}, \theta_{y}) \in D, 
\end{gather*}
where for $m,i \in \mathbb{N}$, $l = 1, \dotsc, h(m,d)$, $j = 1, \dotsc, h(i,d)$, the entries of $\mathfrak{a}_{x,y}(\kappa, \iota)$ are given by
\begin{gather*}
    \mathfrak{a}_{m,i}^{l,j}(\kappa, \iota) = K_{d}^{2}S_{m}^{l}(\theta_{x})\frac{J_{m+\frac{d-2}{2}}(r_{x}\kappa)}{(r_{x}\kappa)^{\frac{d-2}{2}}}S_{i}^{j}(\theta_{y})\frac{J_{i+\frac{d-2}{2}}( r_{y}\iota)}{(r_{y}\iota)^{\frac{d-2}{2}}}.
\end{gather*}
Then, similar to $L^{0}_{D \times D}$, we define $\mathcal{L}^{0}_{D \times D}$ as the space of sequences of functions which are of the form  
\begin{gather*}
    a(\kappa, \iota) = \sum_{p,q}\beta_{pq}\mathfrak{a}_{x_{p},y_{q}}(\kappa, \iota),
\end{gather*}
where $x_{p}, y_{q} \in D$ and $\beta_{pq}$ are real coefficients. Finally, we define $\mathcal{L}_{D \times D}(\Phi_{1} \times \Phi_{2})$ as the closure of $\mathcal{L}^{0}_{D \times D}$ with respect to the inner product 
\begin{gather*}
    \langle a, b \rangle_{\Phi_{1} \times \Phi_{2}} \coloneqq \sum_{m=0}^{\infty}\sum_{l=1}^{h(m,d)}\!\int_{0}^{\infty}\sum_{i=0}^{\infty}\sum_{j=1}^{h(i,d)}\!\int_{0}^{\infty}a_{m,i}^{l,j}(\kappa, \iota)b_{m,i}^{l,j}(\kappa, \iota)\Phi_{1}(d \kappa)\Phi_{2}(d \iota).
\end{gather*}
Using \eqref{coresp1}, we observe that
\begin{equation} \label{aha}
    \langle  \mathfrak{e}_{x_{p},y_{q}}, \mathfrak{e}_{x,y}  \rangle_{F_{1} \times F_{2}} = \langle \mathfrak{a}_{x_{p},y_{q}}, \mathfrak{a}_{x,y}  \rangle_{\Phi_{1} \times \Phi_{2}}.
\end{equation}
According to Theorem~8 in Section~III.3 of \cite{GaussianRandomProcesses}, the Gaussian measures $\mathbb{P}_{1}$ and $\mathbb{P}_{2}$ are equivalent on $\sigma_{D}(\xi)$ if and only if the restriction $\delta$ permits a representation
    \begin{gather*}
    \delta(x,y) = \langle \Psi, \mathfrak{e}_{x,y}  \rangle_{F_{1} \times F_{2}},
    \end{gather*}
with $\Psi \in L_{D \times D}(F_{1} \times F_{2})$. Thus, by definition of $L_{D \times D}(F_{1} \times F_{2})$, we write $\Psi = \lim_{n \to \infty}\Psi_{n}$, where $(\Psi_{n}) \subset L^{0}_{D \times D}$. Using \eqref{aha}, we have that
\begin{align*}
    \delta(x,y) &= \langle\Psi,\mathfrak{e}_{x,y} \rangle_{F_{1} \times F_{2}} \\
    &= \lim_{n \to \infty} \langle \Psi_{n}, \mathfrak{e}_{x,y}  \rangle_{F_{1} \times F_{2}} \\ 
    &= \lim_{n \to \infty} \langle \psi_{n}, \mathfrak{a}_{x,y}  \rangle_{\Phi_{1} \times \Phi_{2}} \\
    &= \langle \psi, \mathfrak{a}_{x,y}  \rangle_{\Phi_{1} \times \Phi_{2}},
\end{align*}
with $(\psi_{n}) \subset \mathcal{L}^{0}_{D \times D}$ such that $\lim_{n \to \infty}\psi_{n} = \psi \in \mathcal{L}_{D \times D}(\Phi_{1} \times \Phi_{2})$. Hence, we have shown that $\mathbb{P}_{1}$ and $\mathbb{P}_{2}$ are equivalent on $\sigma_{D}(\xi)$ if and only if $\delta$ is representable as
\begin{gather*} \label{product}
    \delta(x,y) = K_{d}^{2}\sum_{m=0}^{\infty}\sum_{l=1}^{h(m,d)}\!\int_{0}^{\infty}\sum_{i=0}^{\infty}\sum_{j=1}^{h(i,d)}\!\int_{0}^{\infty}S_{m}^{l}(\theta_{x})\frac{J_{m+\frac{d-2}{2}}(r_{x}\kappa)}{(r_{x}\kappa)^{\frac{d-2}{2}}}S_{i}^{j}(\theta_{y})\frac{J_{i+\frac{d-2}{2}}(r_{y}\iota)}{( r_{y}\iota)^{\frac{d-2}{2}}}\psi_{m,i}^{l,j}(\kappa, \iota)\Phi_{1}(d \kappa)\Phi_{2}(d \iota),
\end{gather*}
where 
\begin{gather*}
\sum_{m=0}^{\infty}\sum_{l=1}^{h(m,d)}\!\int_{0}^{\infty}\sum_{i=0}^{\infty}\sum_{j=1}^{h(i,d)}\!\int_{0}^{\infty}\lvert \psi_{m,i}^{l,j}(\kappa, \iota)\rvert^{2}\Phi_{1}(d \kappa)\Phi_{2}(d \iota) < \infty.
\end{gather*}
Then, we can conclude the proof using the orthonormality property of the spherical harmonics.
\end{proof}

Notice, if Assumptions~\ref{a2} and~\ref{a3} are satisfied, then, for $\ell = 1,2$, since $k_{\ell}$ is the Fourier transform of $f_{\ell}$, and $k_{\ell}$ is assumed to be radial, we must conclude that $f_{\ell}$ is radial itself. We denote its radial version with $g_{\ell}$, i.e., $f_{\ell}(x) = g_{\ell}(\lVert x \rVert)$. The analog to Theorem~\ref{thm1} for isotropic random functions reads as follows:

\begin{theorem} \label{thm2}
Suppose that Assumptions~\ref{a2} and~\ref{a3} are satisfied where $f_{1}$ and $f_{2}$ are bounded on $\Rd$. Then, the Gaussian measures $\mathbb{P}_{1}$ and $\mathbb{P}_{2}$ are equivalent on $\sigma_{D}(\xi)$ if and only if, for any pair $m,i \in \mathbb{N}$, $l = 1, \dotsc, h(m,d)$, $j = 1, \dotsc, h(i,d)$, the scaled difference 
\begin{gather*}
       \mathfrak{d}_{m,i}^{l,j}(r_{1}, r_{2}) \coloneqq r_{1}^{\frac{d-1}{2}}r_{2}^{\frac{d-1}{2}}\delta_{m,i}^{l,j}(r_{1}, r_{2}), \quad r_{1}, r_{2} \in D_{+},
\end{gather*}
satisfies the following properties:
\begin{enumerate}[label=(\alph*)]
    \item \label{thm2.1} There exists extension $\prescript{\scalebox{0.5}{$\uparrow$}}{}{}\mathfrak{d}_{m,i}^{l,j}$ of $\mathfrak{d}_{m,i}^{l,j}$ to $\Rplus \times \Rplus$ which is square-integrable, i.e., $$\int_{0}^{\infty}\int_{0}^{\infty}\lvert \prescript{\scalebox{0.5}{$\uparrow$}}{}{}\mathfrak{d}_{m,i}^{l,j}(r_{1}, r_{2})\rvert^{2}dr_{1}dr_{2} < \infty\, ;$$
    \item \label{thm2.2} The two-dimensional Hankel transform
    \begin{gather*}
        h_{m,i}^{l,j}(\kappa, \iota) = \int_{0}^{\infty}\int_{0}^{\infty}\prescript{\scalebox{0.5}{$\uparrow$}}{}{}\mathfrak{d}_{m,i}^{l,j}(r_{1}, r_{2})\sqrt{r_{1} \kappa}\sqrt{\vphantom{r_{1}\kappa}r_{2} \iota}J_{m+\frac{d-2}{2}}(r_{1} \kappa)J_{i+\frac{d-2}{2}}(r_{2} \iota)dr_{1}dr_{2},
    \end{gather*}
    of $\prescript{\scalebox{0.5}{$\uparrow$}}{}{}\mathfrak{d}_{m,i}^{l,j}$ satisfies 
\begin{gather*}
\sum_{m=0}^{\infty}\sum_{l=1}^{h(m,d)}\!\int_{0}^{\infty}\sum_{i=0}^{\infty}\sum_{j=1}^{h(i,d)}\!\int_{0}^{\infty}\frac{\lvert h_{m,i}^{l,j}(\kappa, \iota)\rvert^{2}}{g_{1}(\kappa) g_{2}(\iota)} d\kappa d\iota < \infty.
\end{gather*}
\end{enumerate}
\end{theorem}

\begin{proof}
We first notice that for any $b \geq 0$, by Assumption~\ref{a2},
\begin{gather*}
\Phi_{\ell}([0,b)) = F_{\ell}(B_{b}(0)) = \int_{\{x \colon \lVert x \rVert < b\}}f_{\ell}(x)dx, \quad \ell = 1,2.
\end{gather*}
Since for $\ell = 1,2$, $f_{\ell}$ is radial, we use spherical coordinates and get
\begin{equation} \label{radialmeasure1}
    \Phi_{\ell}([0,b)) = \frac{2 \pi^{\frac{d}{2}}}{\Gamma\big(\frac{d}{2}\big)}\int_{0}^{b}\kappa^{d-1}g_{\ell}(\kappa)d\kappa, \quad \ell = 1,2.
\end{equation}
Assume that $\mathbb{P}_{1}$ and $\mathbb{P}_{2}$ are equivalent on $\sigma_{D}(\xi)$. We apply Lemma~\ref{intermediateLemma} and conclude that
\begin{gather*}
    \mathfrak{d}_{m,i}^{l,j}(r_{1}, r_{2}) = \frac{2^{d+1}\pi^{\frac{3d}{2}}}{\Gamma\big(\frac{d}{2}\big)}\int_{0}^{\infty}\int_{0}^{\infty}\kappa^{\frac{d-1}{2}}g_{1}(\kappa)\iota^{\frac{d-1}{2}}g_{2}(\iota)\psi_{m,i}^{l,j}(\kappa, \iota)\sqrt{r_{1}\kappa}\sqrt{\vphantom{r_{1}\kappa}r_{2}\iota}J_{m+\frac{d-2}{2}}(\kappa r_{1})J_{i+\frac{d-2}{2}}(\iota r_{2}) d\kappa d \iota,
\end{gather*}
whenever $r_{1}, r_{2} \in D_{+}$, for some $\psi_{m,i}^{l,j}$ which satisfies \eqref{squarepsi}. Notice that since $f_{1}$ and $f_{2}$ are real valued and bounded on $\Rd$, the respective radial versions $g_{1}$ and $g_{2}$ must be real valued and bounded on $\Rplus$. Thus, if we define 
\begin{gather*}
    h_{m,i}^{l,j}(\kappa, \iota) \coloneqq \frac{2^{d+1}\pi^{\frac{3d}{2}}}{\Gamma\big(\frac{d}{2}\big)}\kappa^{\frac{d-1}{2}}g_{1}(\kappa)\iota^{\frac{d-1}{2}}g_{2}(\iota)\psi_{m,i}^{l,j}(\kappa, \iota), \quad \kappa, \iota \in \Rplus,
\end{gather*}
using \eqref{squarepsi}, together with \eqref{radialmeasure1}, we have that
\begin{gather*}
    \int_{0}^{\infty}\int_{0}^{\infty}\lvert h_{m,i}^{l,j}(\kappa, \iota) \rvert^{2} d\kappa d\iota < \infty.
\end{gather*}
Therefore, we define $\prescript{\scalebox{0.5}{$\uparrow$}}{}{}\mathfrak{d}_{m,i}^{l,j}$ on $\Rplus \times \Rplus$, as the two-dimensional Hankel transform of $h_{m,i}^{l,j}(\kappa, \iota)$, which is square-integrable (see Corollary~6.1 in \cite{Hankel}). This then proves \ref{thm2.1} of Theorem~\ref{thm2}. In addition, by \eqref{squarepsi}, also \ref{thm2.2} of Theorem~\ref{thm2} must be satisfied. To prove the other direction, suppose that \ref{thm2.1} and \ref{thm2.2} of Theorem~\ref{thm2} are satisfied. Then by \ref{thm2.1}, since $\prescript{\scalebox{0.5}{$\uparrow$}}{}{}\mathfrak{d}_{m,i}^{l,j}$ is square-integrable on $\Rplus \times \Rplus$, the two-dimensional Hankel transform 
\begin{gather*}
    \tilde{h}_{m,i}^{l,j} \quad \text{of} \quad \bigg(\frac{2^{d+1}\pi^{\frac{3d}{2}}}{\Gamma\big(\frac{d}{2}\big)}\bigg)^{-1}\prescript{\scalebox{0.5}{$\uparrow$}}{}{}\mathfrak{d}_{m,i}^{l,j}
\end{gather*}
exists and is square-integrable on $\Rplus \times \Rplus$. Therefore, on $D_{+} \times D_{+}$, we have that 
\begin{gather*}
    r_{1}^{\frac{d-1}{2}}r_{2}^{\frac{d-1}{2}}\delta_{m,i}^{l,j}(r_{1}, r_{2}) = \frac{2^{d+1}\pi^{\frac{3d}{2}}}{\Gamma\big(\frac{d}{2}\big)}\int_{0}^{\infty}\int_{0}^{\infty}\tilde{h}_{m,i}^{l,j}(\kappa, \iota)\sqrt{r_{1} \kappa}\sqrt{\vphantom{r_{1}\kappa}r_{2} \iota}J_{m+\frac{d-2}{2}}(r_{1} \kappa)J_{i+\frac{d-2}{2}}(r_{2} \iota)d\kappa d\iota.
\end{gather*}
We set
\begin{gather*}
    \psi_{m,i}^{l,j}(\kappa, \iota) \coloneqq \frac{\tilde{h}_{m,i}^{l,j}(\kappa, \iota)}{\kappa^{\frac{d-1}{2}}g_{1}(\kappa)\iota^{\frac{d-1}{2}}g_{2}(\iota)},
\end{gather*}
and obtain
\begin{gather*}
r_{1}^{\frac{d-1}{2}}r_{2}^{\frac{d-1}{2}}\delta_{m,i}^{l,j}(r_{1}, r_{2}) = K_{d}^{2}\int_{0}^{\infty}\int_{0}^{\infty}\sqrt{r_{1}r_{2}}\frac{J_{m+\frac{d-2}{2}}(r_{1}\kappa)}{\kappa^{\frac{d-2}{2}}}\frac{J_{i+\frac{d-2}{2}}( r_{2}\iota)}{\iota^{\frac{d-2}{2}}}\psi_{m,i}^{l,j}(\kappa, \iota)\Phi_{1}(d \kappa)\Phi_{2}(d \iota),
\end{gather*}
with $r_{1}, r_{2} \in D_{+}$. Further, by \ref{thm2.2} of Theorem~\ref{thm2}, we have that
\begin{gather*}
\sum_{m=0}^{\infty}\sum_{l=1}^{h(m,d)}\!\int_{0}^{\infty}\sum_{i=0}^{\infty}\sum_{j=1}^{h(i,d)}\!\int_{0}^{\infty}\lvert \psi_{m,i}^{l,j}(\kappa, \iota) \rvert^{2}\Phi_{1}(d\kappa)\Phi_{2}(d\iota) < \infty.
\end{gather*}
Finally, we can conclude that $\mathbb{P}_{1}$ and $\mathbb{P}_{2}$ are equivalent on $\sigma_{D}(\xi)$ under application of Lemma~\ref{intermediateLemma}.
\end{proof}

Similar to Theorem~\ref{thm1_extension}, the equivalence of $\mathbb{P}_{1}$ and $\mathbb{P}_{2}$ allows for an extension of $\delta_{m,i}^{l,j}$ that is continuous. To arrive there, we introduce the set 
\begin{gather*}
    N \coloneqq \big\{(r_{1},r_{2})\in \Rplus \times \Rplus \colon (r_{1},r_{2}) \in \{0\}\times\Rplus \cup \Rplus \times \{0\} \big\}.
\end{gather*}

\begin{theorem} \label{thm2_extension}
    Suppose that Assumptions~\ref{a2} and~\ref{a3} are satisfied where $f_{1}$ and $f_{2}$ are bounded on $\Rd$. Then, if the Gaussian measures $\mathbb{P}_{1}$ and $\mathbb{P}_{2}$ are equivalent on $\sigma_{D}(\xi)$, for any pair $m,i \in \mathbb{N}$, $l = 1, \dotsc, h(m,d)$, $j = 1, \dotsc, h(i,d)$, there exists a continuous extension $\prescript{\scalebox{0.5}{$\uparrow$}}{}{}\delta_{m,i}^{l,j}$ of $\delta_{m,i}^{l,j}$, from $D_{+} \times D_{+} \setminus N$ to $\Rplus \times \Rplus \setminus N$, which is such that 
    \begin{gather*}
        r_{1}^{\frac{d-1}{2}}r_{2}^{\frac{d-1}{2}}\prescript{\scalebox{0.5}{$\uparrow$}}{}{}\delta_{m,i}^{l,j}(r_{1}, r_{2})
    \end{gather*}
is square-integrable on $\Rplus \times \Rplus$.
\end{theorem}

\begin{proof}
   Let us fix a pair $m,i \in \mathbb{N}$, $l = 1, \dotsc, h(m,d)$, $j = 1, \dotsc, h(i,d)$. From the proof of Theorem~\ref{thm2} we see that the equivalence of $\mathbb{P}_{1}$ and $\mathbb{P}_{2}$ on $\sigma_{D}(\xi)$ implies that
\begin{gather*}
    \delta_{m,i}^{l,j}(r_{1}, r_{2}) = \frac{2^{d+1}\pi^{\frac{3d}{2}}}{\Gamma\big(\frac{d}{2}\big)}\int_{0}^{\infty}\int_{0}^{\infty}\prescript{l,j}{m,i}{w_{r_{1},r_{2}}}(\kappa, \iota) d\kappa d \iota, \quad (r_{1},r_{2}) \in D_{+} \times D_{+} \setminus N,
\end{gather*}
with
\begin{gather*}
    \prescript{l,j}{m,i}{w_{r_{1},r_{2}}}(\kappa, \iota) = \frac{\kappa^{\frac{d-1}{2}}g_{1}(\kappa)\iota^{\frac{d-1}{2}}g_{2}(\iota)}{r_{1}^{\frac{d-1}{2}}r_{2}^{\frac{d-1}{2}}}\psi_{m,i}^{l,j}(\kappa, \iota)\sqrt{r_{1}\kappa}\sqrt{\vphantom{r_{1}\kappa}r_{2}\iota}J_{m+\frac{d-2}{2}}(\kappa r_{1})J_{i+\frac{d-2}{2}}(\iota r_{2}).
\end{gather*}
We recall that $\psi_{m,i}^{l,j}$ is square-integrable with respect to
\begin{gather*}
    \Phi_{1}(d \kappa)\Phi_{2}(d \iota) = \bigg(\frac{2\pi^{\frac{d}{2}}}{\Gamma\big(\frac{d}{2}\big)}\bigg)^{2}\kappa^{d-1}g_{1}(\kappa)\iota^{d-1}g_{2}(\iota)d\kappa d\iota.
\end{gather*}
As with the proof of Theorem~\ref{thm1_extension}, since $\Phi_{1}$ and $\Phi_{2}$ are bounded, we conclude that 
\begin{equation} \label{absintegrable}
    \int_{0}^{\infty}\int_{0}^{\infty}\lvert \psi_{m,i}^{l,j}(\kappa, \iota)\rvert\Phi_{1}(d \kappa)\Phi_{2}(d \iota) < \infty,
\end{equation}
as well. Then, for any pair $(r_{1}, r_{2}) \in \Rplus\times \Rplus \setminus N$, we define, for any $\kappa, \iota \in \Rplus$,
\begin{gather*}
    \prescript{l,j}{m,i}{W_{r_{1},r_{2}}}(\kappa, \iota) \coloneqq \frac{\kappa^{\frac{d-1}{2}}g_{1}(\kappa)\iota^{\frac{d-1}{2}}g_{2}(\iota)}{r_{1}^{\frac{d-1}{2}}r_{2}^{\frac{d-1}{2}}}\psi_{m,i}^{l,j}(\kappa, \iota)\sqrt{r_{1}\kappa}\sqrt{\vphantom{r_{1}\kappa}r_{2}\iota}J_{m+\frac{d-2}{2}}(\kappa r_{1})J_{i+\frac{d-2}{2}}(\iota r_{2}).
\end{gather*}
We remark that for fixed $\kappa, \iota \in \Rplus$, $(r_{1}, r_{2}) \mapsto \prescript{l,j}{m,i}{W_{r_{1},r_{2}}}(\kappa, \iota)$ is continuous on $\Rplus\times \Rplus \setminus N$. Then, using Lommel's expression for the Bessel function of the first kind (see Section 3.3 of \cite{Watson}), we estimate
\begin{align*}
    J_{m+\frac{d-2}{2}}(\kappa r_{1}) &= \frac{\big(\frac{1}{2}\kappa r_{1}\big)^{m+\frac{d-2}{2}}}{\Gamma\big(m+\frac{d-2}{2} + \frac{1}{2}\big)\Gamma\big(\frac{1}{2}\big)}\int_{0}^{\pi}\operatorname{cos}(\kappa r_{1} \operatorname{cos}(\theta))\operatorname{sin}^{2\big(m+\frac{d-2}{2}\big)}(\theta)d\theta \\
    &\leq \frac{\big(\frac{1}{2}\kappa r_{1}\big)^{\frac{d-1}{2}}\Gamma\big(m\big)}{\Gamma\big(m+\frac{d-2}{2} + \frac{1}{2}\big)}\frac{\big(\frac{1}{2}\kappa r_{1}\big)^{m-\frac{1}{2}}}{\Gamma\big(m\big)\Gamma\big(\frac{1}{2}\big)}\int_{0}^{\pi}\operatorname{cos}(\kappa r_{1} \operatorname{cos}(\theta))\operatorname{sin}^{2\big(m-\frac{1}{2}\big)}(\theta)d\theta\\
    &= \frac{\big(\frac{1}{2}\kappa r_{1}\big)^{\frac{d-1}{2}}\Gamma\big(m\big)}{\Gamma\big(m+\frac{d-2}{2} + \frac{1}{2}\big)}J_{m-\frac{1}{2}}(\kappa r_{1}).
\end{align*}
Hence, if one applies a similar estimate to $J_{i+(d-2)/2}(\iota r_{2})$, we obtain that
\begin{gather*}
    \prescript{l,j}{m,i}{W_{r_{1},r_{2}}}(\kappa, \iota) \leq C\kappa^{d-1}g_{1}(\kappa)\iota^{d-1}g_{2}(\iota)\psi_{m,i}^{l,j}(\kappa, \iota)\sqrt{r_{1}\kappa}\sqrt{\vphantom{r_{1}\kappa}r_{2}\iota}J_{m-\frac{1}{2}}(\kappa r_{1})J_{i-\frac{1}{2}}(\iota r_{2}),
\end{gather*}
where $C$ is some fixed constant, independent of $r_{1},r_{2}$ and $\kappa, \iota$. In addition, since the functions $z \mapsto \sqrt{z}J_{m-1/2}(z)$ and $z \mapsto \sqrt{z}J_{i-1/2}(z)$ are bounded on $(0, \infty)$ (see \cite{Hankel1}), we get for any $(r_{1}, r_{2}) \in \Rplus\times \Rplus \setminus N$, 
\begin{gather*}
    \frac{2^{d+1}\pi^{\frac{3d}{2}}}{\Gamma\big(\frac{d}{2}\big)}\lvert \prescript{l,j}{m,i}{W_{r_{1},r_{2}}}(\kappa, \iota) \rvert \leq \widetilde{C}\frac{2^{d+1}\pi^{\frac{3d}{2}}}{\Gamma\big(\frac{d}{2}\big)}\kappa^{d-1}g_{1}(\kappa)\iota^{d-1}g_{2}(\iota)\psi_{m,i}^{l,j}(\kappa, \iota),
\end{gather*}
where $\widetilde{C}$ is independent of $r_{1},r_{2}$ and $\kappa, \iota$. The function on the right-hand side of the latter equation is absolutely integrable because of \eqref{absintegrable}. Hence, we set
\begin{gather*}
    \prescript{\scalebox{0.5}{$\uparrow$}}{}{}\delta_{m,i}^{l,j}(r_{1}, r_{2}) \coloneqq \frac{2^{d+1}\pi^{\frac{3d}{2}}}{\Gamma\big(\frac{d}{2}\big)}\int_{0}^{\infty}\int_{0}^{\infty}\prescript{l,j}{m,i}{W_{r_{1},r_{2}}}(\kappa, \iota)d\kappa d\iota, \quad (r_{1},r_{2}) \in \Rplus \times \Rplus \setminus N,
\end{gather*}
and obtain an extension of $\delta_{m,i}^{l,j}$ from $D_{+} \times D_{+} \setminus N$ to $\Rplus \times \Rplus \setminus N$, which is continuous by Lebesgue's dominated convergence theorem. Then, using the same reasoning as in the proof of Theorem~\ref{thm2}, since $f_{1}$ and $f_{2}$ are assumed to be bounded, we extend 
\begin{gather} \label{scaled_extension}
    r_{1}^{\frac{d-1}{2}}r_{2}^{\frac{d-1}{2}}\prescript{\scalebox{0.5}{$\uparrow$}}{}{}\delta_{m,i}^{l,j}(r_{1}, r_{2})
\end{gather}
to the entire $\Rplus \times \Rplus$ by means of the two-dimensional Hankel transform of a square-integrable function. Finally, since $N$ has Lebesgue measure zero, \eqref{scaled_extension} is square-integrable on $\Rplus \times \Rplus$.
\end{proof}

\subsection{Sufficiently dense sampling} \label{sec:dense} 

Upon Theorem~\ref{thm1_extension}, the next result relates the uniform continuity of the covariance functions $c_{1}$ and $c_{2}$ with the orthogonality of the Gaussian measures $\mathbb{P}_{1}$ and $\mathbb{P}_{2}$ when $D$ is dense in $\Rd$.

\begin{theorem} \label{thm3}
    Suppose that Assumption~\ref{a2} is satisfied where $f_{1}$ and $f_{2}$ are bounded on $\Rd$ and such that the set $\{\lambda \in \Rd \colon f_{1}(\lambda) \neq f_{2}(\lambda)\}$ has positive Lebesgue measure. Then, if $c_{1}$ and $c_{2}$ are uniformly continuous on $\Rd \times \Rd$ and $D$ is dense in $\Rd$, the Gaussian measures $\mathbb{P}_{1}$ and $\mathbb{P}_{2}$ are orthogonal on $\sigma_{D}(\xi)$.   
\end{theorem}

\begin{proof}
We recall that on $\Rd \times \Rd$, the difference $c_{1}-c_{2}$ is given by
\begin{equation} \label{entiredifference}
c_{1}(x,y) - c_{2}(x,y) = \int_{\Rd}\euler^{\iu \langle \lambda, x-y\rangle}\big(f_{1}(\lambda)-f_{2}(\lambda)\big)d\lambda.
\end{equation}
Since $\{\lambda \in \Rd \colon f_{1}(\lambda) \neq f_{2}(\lambda)\}$ has positive Lebesgue measure, \eqref{entiredifference} can not be square-integrable on $\Rd \times \Rd$. To see it, we recall that $\xi$ is stationary and observe that
\begin{gather*}
    \int_{\Rd}\lvert c_{1}(x,y) - c_{2}(x,y) \rvert^{2} dx = \int_{\Rd}\lvert k_{1}(x-y) - k_{2}(x-y) \rvert^{2} dx.
\end{gather*}
Then, as the later integral is constant in $y$, we must have that
\begin{gather*}
    \int_{\Rd}\int_{\Rd}\lvert k_{1}(x-y) - k_{2}(x-y) \rvert^{2} dx dy = \infty,
\end{gather*}
unless the $L^{2}$ norm of the difference $k_{1} - k_{2}$ is zero. But, since the Fourier transform is an isometry on $L^{2}(\Rd)$, and $f_{1}-f_{2}$ is assumed to have non-zero $L^{2}$ norm, this case is not possible. Thus, \eqref{entiredifference} is not square-integrable. Still, by assumption, it is continuous on $\Rd \times \Rd$. Furthermore, since $\delta$ is assumed to be uniformly continuous on $D \times D$ and $D$ is dense in $\Rd$, any continuous extension of $\delta$ to $\Rd \times \Rd$ must be given by \eqref{entiredifference}. This concludes the proof under application of Theorem~\ref{thm1_extension}.
\end{proof}

If $\xi$ is also isotropic, the density of $D_{+}$ in $\Rplus$ is sufficient to recover the orthogonality of the Gaussian measures $\mathbb{P}_{1}$ and $\mathbb{P}_{2}$ on $\sigma_{D}(\xi)$. 

\begin{theorem} \label{thm4}
    Suppose that Assumptions~\ref{a2} and~\ref{a3} are satisfied where $f_{1}$ and $f_{2}$ are bounded on $\Rd$ and such that the set $\{\lambda \in \Rd \colon f_{1}(\lambda) \neq f_{2}(\lambda)\}$ has positive Lebesgue measure. Then, if $c_{1}$ and $c_{2}$ are uniformly continuous on $\Rd \times \Rd$ and $D_{+}$ is dense in $\Rplus$, the Gaussian measures $\mathbb{P}_{1}$ and $\mathbb{P}_{2}$ are orthogonal on $\sigma_{D}(\xi)$.
\end{theorem}

\begin{proof}
First of all, using \eqref{radial1}, we see that
\begin{align*}
    \prescript{\scalebox{0.5}{$\uparrow$}}{*}{}\delta_{m}(r_{x}, r_{y}) &\coloneqq \int_{\mathbb{S}^{d-1}}\int_{\mathbb{S}^{d-1}}\big(c_{1}(x,y)-c_{2}(x,y)\big)S_{m}^{l}(\theta_{x})S_{m}^{l}(\theta_{y})d\theta_{x}d\theta_{y} \\
               &= (2\pi)^{d}\int_{0}^{\infty}\frac{J_{m+\frac{d-2}{2}}(\kappa r_{x})}{(\kappa r_{x})^{\frac{d-2}{2}}}\frac{J_{m+\frac{d-2}{2}}(\kappa r_{y})}{(\kappa r_{y})^{\frac{d-2}{2}}} \kappa^{d-1}\big(g_{1}(\kappa) - g_{2}(\kappa)\big)d\kappa,
\end{align*}
is one extension of $\delta_{m,i}^{l,j}$ to $\Rplus \times \Rplus$ for the pair $m = i$. Since $c_{1}$ and $c_{2}$ are uniformly continuous on $\Rd \times \Rd$, $\prescript{\scalebox{0.5}{$\uparrow$}}{*}{}\delta_{m}$ must be uniformly continuous on $\Rplus \times \Rplus$. This follows from the fact that 
\begin{gather*}
    \int_{\mathbb{S}^{d-1}}\lvert S_{m}^{l}(\theta_{x}) \rvert d\theta_{x} \leq \bigg(\frac{h(m,d)2\pi^{\frac{d}{2}}}{\Gamma\big(\frac{d}{2}\big)}\bigg)^{\frac{1}{2}}.
\end{gather*}
See for instance (b) of Corollary~2.9 in \cite{WeissStein}. In particular, $\prescript{\scalebox{0.5}{$\uparrow$}}{*}{}\delta_{m}$ is a continuous extension of $\delta_{m,m}^{l,l}$ from $D_{+} \times D_{+} \setminus N$ to $\Rplus \times \Rplus \setminus N$. But, because of the assumption that $\{\lambda \in \Rd \colon f_{1}(\lambda) \neq f_{2}(\lambda)\}$ has positive Lebesgue measure, it can not be the case that
\begin{equation*} \label{square_radial}
    Q_{m} \coloneqq \int_{0}^{\infty}\int_{0}^{\infty}r_{x}^{d-1}r_{y}^{d-1}\lvert\prescript{\scalebox{0.5}{$\uparrow$}}{*}{}\delta_{m}(r_{x}, r_{y})\rvert^{2}dr_{x}dr_{y} < \infty.
\end{equation*}
To see this, we use the identity
\begin{gather*}
    c_{1}(x,y) - c_{2}(x,y) = \sum_{m=0}^{\infty}\sum_{l=1}^{h(m,d)}S_{m}^{l}(\theta_{x})S_{m}^{l}(\theta_{y})\prescript{\scalebox{0.5}{$\uparrow$}}{*}{}\delta_{m}(r_{x}, r_{y}), \quad x,y \in \Rd,
\end{gather*}
and note that
\begin{align*}
    Q_{m} &= \int_{0}^{\infty}\int_{\mathbb{S}^{d-1}}\int_{0}^{\infty}\int_{\mathbb{S}^{d-1}}\big\lvert c_{1}(x,y) - c_{2}(x,y) \big \rvert^{2} r_{x}^{d-1} r_{y}^{d-1} d\theta_{x} d\theta_{y} d_{r_{x}}d_{r_{y}}\\
          &= \int_{\Rd}\int_{\Rd}\big\lvert c_{1}(x,y) - c_{2}(x,y) \big \rvert^{2} dx dy.
\end{align*}
Then, using the same reasoning as in the proof of Theorem~\ref{thm3}, the latter integral is not finite. Therefore, we found a pair $m=i$ for which $\prescript{\scalebox{0.5}{$\uparrow$}}{*}{}\delta_{m}$ is a continuous extension of $\delta_{m,m}^{l,l}$ from $D_{+} \times D_{+} \setminus N$ to $\Rplus \times \Rplus \setminus N$ which is such that $r_{x}^{(d-1)/2}r_{y}^{(d-1)/2}\prescript{\scalebox{0.5}{$\uparrow$}}{*}{}\delta_{m}(r_{x}, r_{y})$ is not square-integrable on $\Rplus \times \Rplus$. Still, $\delta_{m,m}^{l,l}$ is uniformly continuous on $D_{+} \times D_{+}$. In particular it is uniformly continuous on $D_{+} \times D_{+} \setminus N$. Since $D_{+} \times D_{+}$ is dense in $\Rplus \times \Rplus$, $D_{+} \times D_{+} \setminus N$ is dense in $\Rplus \times \Rplus \setminus N$. Hence, any continuous extension of $\delta_{m,m}^{l,l}$ from $D_{+} \times D_{+} \setminus N$ to $\Rplus \times \Rplus \setminus N$ must be given by $\prescript{\scalebox{0.5}{$\uparrow$}}{*}{}\delta_{m}$. This concludes the proof under application of Theorem~\ref{thm2_extension}.
\end{proof}

From the proofs of Theorems~\ref{thm3} and~\ref{thm4}, it becomes obvious that the assumption that $\{\lambda \in \Rd \colon f_{1}(\lambda) \neq f_{2}(\lambda)\}$ has positive Lebesgue measure can be replaced with the assumption that $\{x \in \Rd \colon k_{1}(x) \neq k_{2}(x)\}$ has positive Lebesgue measure. Notice also that if $D$ is dense in $\Rd$, then clearly $D_{+}$ is dense in $\Rplus$. Of course, the converse is not true. We will see a particular example in the next section. We further remark that if $\xi$ is measurable with respect to some larger $\sigma$-algebra $\mathcal{G}$, i.e., $\mathcal{U} \subset \mathcal{G}$, then Theorems~\ref{thm3} and \ref{thm4} give sufficient conditions to deduce the orthogonality of $\mathbb{P}_{1}$ and $\mathbb{P}_{2}$ on $\{\xi^{-1}(G) \colon G \in \mathcal{G}\}$. This is because the orthogonality of $\mathbb{P}_{1}$ and $\mathbb{P}_{2}$ on $\{\xi^{-1}(G) \colon G \in \mathcal{G}\}$ follows from the orthogonality of $\mathbb{P}_{1}$ and $\mathbb{P}_{2}$ on $\sigma_{D}(\xi)$. To finish this section, we have a closer look at a well-known family of covariance functions.

\begin{example}[Exponential family of covariance functions] \label{exponentialfam}
The exponential family is defined by (see Example~1 on p.\ 115 of \cite{Yaglom})
\begin{gather} \label{exponential}
    \phi_{\theta}(\tau) \coloneqq \sigma^{2}\euler^{-\alpha \tau}, \quad \tau \in \Rplus,
\end{gather}
with $\theta = (\sigma^{2}, \alpha) \in (0,\infty)^{2}$. Clearly, the derivative of \eqref{exponential} with respect to $\tau$ is bounded. In particular, for any $\theta \in (0,\infty)^{2}$, $\tau \mapsto \phi_{\theta}(\tau)$ is Lipschitz continuous. Using the fact that the composition of Lipschitz continuous functions is again Lipschitz continuous, we conclude that for any $\theta \in (0,\infty)^{2}$, $c_{\theta}(x,y) \coloneqq \phi_{\theta}(\lVert x-y \rVert)$, $x,y \in \Rd$, is Lipschitz continuous and thus uniformly continuous on $\Rd \times \Rd$.
The family of radial spectral densities associated with \eqref{exponential} is given by (see \cite{Yaglom}), 
\begin{gather*} \label{spectral_exp}
    g_{\theta}(\kappa) = \sigma^{2}\pi^{-1}\frac{\alpha}{\alpha^{2} + \kappa^{2}}, \quad \kappa \in \Rplus.
\end{gather*}
Hence, for any $\theta \in (0,\infty)^{2}$, $f_{\theta}(x) \coloneqq g_{\theta}(\lVert x \lVert)$, $x \in \Rd$, is bounded. Take $\theta_{1}, \theta_{2} \in (0,\infty)^{2}$ such that $\theta_{1} \neq \theta_{2}$. It follows that $f_{\theta_{1}} \neq f_{\theta_{2}}$ on a set of positive Lebesgue measure. To see it, we can work with \eqref{exponential} and use the fact that the Fourier transform in an isometry $L^{2}(\Rd)$. Explicitly, if $\sigma^{2}_{1} \neq \sigma^{2}_{2}$ and either $\alpha_{1} = \alpha_{2}$ or $\alpha_{1} \neq \alpha_{2}$, we observe that $\phi_{\theta_{1}}(0) \neq \phi_{\theta_{2}}(0)$. But, $\phi_{\theta}$ is continuous at zero, thus, there exists an interval $[0,b)$, $b > 0$, on which $\phi_{\theta_{1}} \neq \phi_{\theta_{2}}$. In the remaining case, i.e., if $\alpha_{1} \neq \alpha_{2}$ but $\sigma^{2}_{1} = \sigma^{2}_{2}$, we can write $\alpha_{1} = \alpha_{2} + s$, with $s \neq 0$. Thus, in this case, we get
\begin{gather*}
    \phi_{\theta_{1}}(\tau) - \phi_{\theta_{2}}(\tau) = \sigma^{2}_{1}\euler^{-\alpha_{2}\tau}\big(\euler^{-s\tau}-1\big),
\end{gather*}
which is non-zero on $\Rplus$. Hence, for any $\theta_{1}, \theta_{2} \in (0,\infty)^{2}$ with $\theta_{1} \neq \theta_{2}$ we have that $\phi_{\theta_{1}} \neq \phi_{\theta_{2}}$ on an interval $[0,b)$, $b > 0$. Therefore, with $k_{\theta}(x) \coloneqq \phi_{\theta}(\lVert x \rVert)$, $x \in \Rd$, $k_{\theta_{1}} \neq k_{\theta_{2}}$ on $B_{b}(0)$. Clearly, $k_{\theta}$ is an element of $L^{2}(\Rd)$. Further, $k_{\theta}$ is the Fourier transform of $f_{\theta}$. Thus, for $\theta_{1} \neq \theta_{2}$, since $k_{\theta_{1}} \neq k_{\theta_{2}}$ on $B_{b}(0)$, we must conclude that the $L^{2}$ norm of $f_{\theta_{1}} - f_{\theta_{2}}$ on $\Rd$ is non-zero. Which implies that $f_{\theta_{1}} \neq f_{\theta_{2}}$ on a set of positive Lebesgue measure. In conclusion, two zero-mean Gaussian measures $\mathbb{P}_{\theta_{1}}$ and $\mathbb{P}_{\theta_{2}}$, $\theta_{1}, \theta_{2} \in (0,\infty)^{2}$, $\theta_{1} \neq \theta_{2}$, with exponential covariance functions $c_{\theta_{1}}$ and $c_{\theta_{2}}$, respectively, are orthogonal on $\sigma_{D}(\xi)$ if $D_{+}$ is dense in $\Rplus$.   
\end{example}

\section{Stochastic sampling} \label{sec:RandomSamplingPaths} 

Let $\{X_{t} \colon t \in T\}$, $T \subset \R$, be a stochastic process defined on a probability space $(\Omega_{\scalebox{0.5}{$X$}}, \mathcal{F}_{\scalebox{0.5}{$X$}}, \mathbb{P}_{\scalebox{0.5}{$X$}})$, taking values in $\Rd$. That is, we consider a random function $\omega_{\scalebox{0.5}{$X$}} \mapsto X(\omega_{\scalebox{0.5}{$X$}})$ with $\Rd$ valued sample paths defined on $T$. We assume that $(\Omega_{\scalebox{0.5}{$X$}}, \mathcal{F}_{\scalebox{0.5}{$X$}}, \mathbb{P}_{\scalebox{0.5}{$X$}})$ is complete in the measure theoretic sense. Further, $X$ starts from the origin, i.e., $X_{t_{0}}(\omega_{\scalebox{0.5}{$X$}}) = 0$, $\omega_{\scalebox{0.5}{$X$}} \in \Omega_{\scalebox{0.5}{$X$}}$, and has continuous sample paths. For a given $\omega_{\scalebox{0.5}{$X$}} \in \Omega_{\scalebox{0.5}{$X$}}$, $X[T](\omega_{\scalebox{0.5}{$X$}})$ denotes the image of $X(\omega_{\scalebox{0.5}{$X$}})$, i.e., $x \in X[T](\omega_{\scalebox{0.5}{$X$}})$ if and only if $x = X_{t}(\omega_{\scalebox{0.5}{$X$}})$ for some $t \in T$. For now, we assume that $\xi$ introduced in Section~\ref{sec:head} is observed along sample paths restricted to $X[T](\omega_{\scalebox{0.5}{$X$}})$, $\omega_{\scalebox{0.5}{$X$}} \in \Omega_{\scalebox{0.5}{$X$}}$. Explicitly, we consider two Gaussian measures $\mathbb{P}_{1}$ and $\mathbb{P}_{2}$ on $\sigma_{X[T](\omega_{\scalebox{0.5}{$X$}})}(\xi)$ which differ only in their covariance functions $c_{1}$ and $c_{2}$. To adapt the notation of the previous sections, we put $X[T]_{+}(\omega_{\scalebox{0.5}{$X$}}) \coloneqq \{\lVert x \rVert \colon x \in X[T](\omega_{\scalebox{0.5}{$X$}})\}$, $\omega_{\scalebox{0.5}{$X$}} \in \Omega_{\scalebox{0.5}{$X$}}$. Following Theorem~\ref{thm3} we obtain: 

\begin{corollary} \label{c0}
    Suppose that the assumptions of Theorem~\ref{thm3} are satisfied with $c_{1}$ and $c_{2}$ uniformly continuous on $\Rd \times \Rd$. Then, if $\mathbb{P}_{\scalebox{0.5}{$X$}}(X[T] \text{ is dense in } \Rd) = 1$, we have that 
    \begin{equation} \label{asortho}
    \mathbb{P}_{\scalebox{0.5}{$X$}}\big(\mathbb{P}_{1} \perp \mathbb{P}_{2} \text{ on $\sigma_{X[T]}(\xi)$}\big) = 1,
    \end{equation}
    i.e., the Gaussian measures $\mathbb{P}_{1}$ and $\mathbb{P}_{2}$ are orthogonal on $\sigma_{X[T]}(\xi)$ $\mathbb{P}_{\scalebox{0.5}{$X$}}$ a.s. 
\end{corollary}

We note that the set $\{X[T] \text{ is dense in } \Rd\}$ is a member of $\mathcal{F}_{\scalebox{0.5}{$X$}}$ since $X$ takes values in $\Rd$ and has continuous sample paths on $T$. Further, since $(\Omega_{\scalebox{0.5}{$X$}}, \mathcal{F}_{\scalebox{0.5}{$X$}}, \mathbb{P}_{\scalebox{0.5}{$X$}})$ is assumed to be complete, under the assumptions of Corollary~\ref{c0}, $\{\mathbb{P}_{1} \perp \mathbb{P}_{2} \text{ on $\sigma_{X[T]}(\xi)$}\}$ is a member of $\mathcal{F}_{\scalebox{0.5}{$X$}}$ as well. Similarly, for the isotropic case, the next result is deduced from Theorem~\ref{thm4}.

\begin{corollary} \label{c0.1}
    Suppose that the assumptions of Theorem~\ref{thm4} are satisfied with $c_{1}$ and $c_{2}$ uniformly continuous on $\Rd \times \Rd$. Then, if $\mathbb{P}_{\scalebox{0.5}{$X$}}(X[T]_{+} \text{ is dense in } \Rplus) = 1$, \eqref{asortho} is satisfied.
\end{corollary}

As we have assumed that $X$ starts from the origin, $X[T]_{+}$ is actually equal to $\Rplus$, with $\mathbb{P}_{\scalebox{0.5}{$X$}}$ probability one, if $X$ has sample paths that are almost surely unbounded. This is summarized in the following lemma:

\begin{lemma} \label{path_connected}
    Suppose that $X[T]$ is unbounded $\mathbb{P}_{\scalebox{0.5}{$X$}}$ a.s. Then $\mathbb{P}_{\scalebox{0.5}{$X$}}(X[T]_{+} = \Rplus) = 1$.  
\end{lemma}

\begin{proof}
    Let $\omega_{\scalebox{0.5}{$X$}} \in \{X[T]\text{ is unbounded}\}$. Clearly $X[T]_{+}(\omega_{\scalebox{0.5}{$X$}}) \subset \Rplus$. For the other direction, since the sample paths of $X$ are continuous, we have that for any $\omega_{\scalebox{0.5}{$X$}} \in \Omega_{\scalebox{0.5}{$X$}}$, $X[T](\omega_{\scalebox{0.5}{$X$}})$ is path-connected. Consider any $r \in \Rplus$ and the neighborhood $B_{r}(0)$ around the origin. Since $\omega_{\scalebox{0.5}{$X$}} \in \{X[T]\text{ is unbounded}\}$, we must conclude that there exists $v$ such that $v \in X[T](\omega_{\scalebox{0.5}{$X$}}) \setminus B_{r}(0)$. But since $X[T](\omega_{\scalebox{0.5}{$X$}})$ is path-connected, there also exists $v^{\prime}$ such that $v^{\prime} \in \partial B_{r}(0) \cap X[T](\omega_{\scalebox{0.5}{$X$}})$. Therefore we have that $\lVert v^{\prime} \rVert = r$, which shows that $\Rplus \subset X[T]_{+}(\omega_{\scalebox{0.5}{$X$}})$.
\end{proof}

Using Lemma~\ref{path_connected}, we have proven the following theorem:

\begin{theorem} \label{thm5}
    Suppose that the assumptions of Theorem~\ref{thm4} are satisfied with $c_{1}$ and $c_{2}$ uniformly continuous on $\Rd \times \Rd$. Then, if $X[T]$ is unbounded $\mathbb{P}_{\scalebox{0.5}{$X$}}$ a.s., we have
    \begin{gather*}
        \mathbb{P}_{\scalebox{0.5}{$X$}}\big(\mathbb{P}_{1} \perp \mathbb{P}_{2} \text{ on $\sigma_{X[T]}(\xi)$}\big) = 1.
    \end{gather*}
\end{theorem}

\begin{example}[Gaussian random fields sampled along Brownian paths] \label{exampleBM}
It is well known that if we let $X = B$, with $T = \R_{+}$, a $d$-dimensional Brownian motion starting from the origin, then for $d \geq 2$, the Lebesgue measure of $B[\R_{+}]$ is zero with $\mathbb{P}_{\scalebox{0.5}{$B$}}$ probability one. This is shown in Th\'{e}or\`{e}me 53, p.\ 240, of \cite{levy} for the case where $d = 2$. A more recent proof, for the general case ($d \geq 2$), is given in the second paragraph of p.\ 197 in \cite{LeGall}. Still, it is true that for $d =1,2$, $B[\R_{+}]$ is dense in $\R$, $\R^{2}$, respectively, with $\mathbb{P}_{\scalebox{0.5}{$B$}}$ probability one (see Propositions~2.14 and~7.16 in \cite{LeGall}). For the case where $d \geq 3$, we know that $\mathbb{P}_{\scalebox{0.5}{$B$}}$ a.s.\ $\lim_{t \to \infty} \lVert B_{t} \rVert = \infty$ (see Theorem~7.17 in \cite{LeGall}). Further, for any $d \geq 1$, the sample paths of $B$ are continuous (see Definitions 2.12 and 2.24 of \citep{LeGall}). In conclusion, given any $d\geq1$, under sufficient conditions on $c_{1}$ and $c_{2}$ (see Corollary~\ref{c0} and Theorem~\ref{thm5}), with $\mathbb{P}_{\scalebox{0.5}{$B$}}$ probability one, Gaussian measures $\mathbb{P}_{1}$ and $\mathbb{P}_{2}$ are orthogonal on $\sigma_{B[\R_{+}]}(\xi)$, i.e., when $\xi$ is sampled along the paths of $B$.
\end{example}

\section{Inference on random fields} \label{sec:MLestimation}

In this section, we let $D = \{x_{i} \colon i \in \natnum\}$ be a fixed sequence of coordinates in $\mathbb{R}^{d}$. That is, we consider 
\begin{gather*}
    \sigma_{D}(\xi) = \big\{\xi^{-1}(U) \colon U \in \mathcal{U}\big\},
\end{gather*}
where now $\mathcal{U} = \sigma\big(\cup_{n=1}^{\infty}\mathcal{U}_{n}\big)$, with
\begin{gather*}
    \mathcal{U}_{n} = \sigma\big(\{C_{x_{1}, \dotsc, x_{n}}(B_{n}) \colon B_{n} \in \mathfrak{B}(\Rn)\}\big).
\end{gather*}
Let $Y_{n} \coloneqq (\xi_{x_{1}}, \dotsc, \xi_{x_{n}})$, $n \in \natnum$. Then, 
\begin{gather*}
    \sigma_{D}(\xi) = \sigma\big(\cup_{n=1}^{\infty}\sigma(Y_{n})\big), \quad \sigma(Y_{n}) = \{Y_{n}^{-1}(B_{n}) \colon B_{n} \in \mathfrak{B}(\Rn)\}.
\end{gather*}
Let $\Theta \subset \mathbb{R}^{p}$. Suppose that $\mathbb{P}_{\theta}$, $\theta \in \Theta$, is a family of Gaussian measures defined on $\sigma_{D}(\xi)$. We remain in the setting of Section~\ref{sec:EquivOrtho}, i.e., for any two $\theta_{1}, \theta_{2} \in \Theta$, under $\mathbb{P}_{\theta_{1}}$ and $\mathbb{P}_{\theta_{2}}$, $\xi$ is stationary and items \ref{a1.1} and \ref{a1.2} are satisfied. Suppose that there exists $\theta_{0} \in \Theta$ such that the true distribution of $\xi$ is obtained from $\mathbb{P}_{\theta_{0}}$. It is further assumed that there exists a neighborhood of $\theta_{0}$ in $\Theta$. In the framework of parameter estimation, $\theta_{0}$ is treated as the unknown and $\Theta$ is regarded as the parameter space. A maximum likelihood (ML) estimator for $\theta_{0}$ is defined to be any sequence of random variables $(\hat{\theta}_{n})$, which is such that for any $n \in \natnum$,
\begin{gather*}
    \hat{\theta}_{n}(\omega) \in \argmax_{\theta \in \Theta}p_{n}(\theta)(\omega), \quad \omega \in \Omega,
\end{gather*}
where 
\begin{gather*}
    p_{n}(\theta) = \frac{1}{\sqrt{(2\pi)^{n}\det\Sigma_{n}(\theta)}}\euler^{-\frac{1}{2}Y_{n}^{\mathrm{t}}\Sigma_{n}(\theta)^{-1}Y_{n}},
\end{gather*}
with
\begin{gather*}
    \Sigma_{n}(\theta) = \big[c_{\theta}(x_{i}, x_{j})\big]_{1 \leq i,j \leq n}.
\end{gather*}
The function $\theta \mapsto p_{n}(\theta)(\omega)$ is called the likelihood function, the probability density function of $Y_{n}$, regarded as a function of $\theta$. The sequence $(\hat{\theta}_{n})$ is said to be strongly consistent for $\theta_{0}$ if 
\begin{gather*}
    \mathbb{P}_{\theta_{0}}\Big(\hat{\theta}_{n} \xrightarrow[]{ n \to \infty } \theta_{0}\Big) = 1.
\end{gather*}
We say that the family $\mathbb{P}_{\theta}$, $\theta \in \Theta$, is a family of orthogonal Gaussian measures on $\sigma_{D}(\xi)$ if for any two $\theta_{1}, \theta_{2} \in \Theta$, with $\theta_{1} \neq \theta_{2}$, $\mathbb{P}_{\theta_{1}}$ and $\mathbb{P}_{\theta_{2}}$ are orthogonal on $\sigma_{D}(\xi)$. The next result is inspired by \cite{Zhang} (see the proof of Theorem~3).
\begin{theorem} \label{consistency}
Let $\Theta$ be closed and convex. Assume that a ML estimator for $\theta_{0}$ exists and
\begin{gather*}
        \theta \mapsto \varphi_{n}(\theta)(\omega) \coloneqq \frac{p_{n}(\theta)(\omega)}{p_{n}(\theta_{0})(\omega)}, \quad \omega \in \Omega,
\end{gather*}
is continuous on $\Theta$. Suppose that there exists $N \in \natnum$ such that for any $n \geq N$ and $\omega \in \Omega$, $\theta \mapsto \varphi_{n}(\theta)(\omega)$ is log-concave on $\Theta$. Then, if $\mathbb{P}_{\theta}$, $\theta \in \Theta$, is a family of orthogonal Gaussian measures on $\sigma_{D}(\xi)$, $(\hat{\theta}_{n})$ is strongly consistent. 
\end{theorem}

\begin{proof}
First of all, since $\mathbb{P}_{\theta}$, $\theta \in \Theta$, is a family of orthogonal Gaussian measures on $\sigma_{D}(\xi)$, we have that with $\mathbb{P}_{\theta_{0}}$ probability one, $(\varphi_{n}(\theta))$ converges pointwise to zero whenever $\theta \neq \theta_{0}$, i.e., 
\begin{gather*}
    \mathbb{P}_{\theta_{0}}\Big(\varphi_{n}(\theta) \xrightarrow[]{n \to \infty} 0\Big) = 1, \quad \theta \in \Theta, \ \theta \neq \theta_{0}.
\end{gather*}
This follows from the fact that the sequence $(\varphi_{n}(\theta))$ forms a martingale on $(\Omega, \sigma_{D}(\xi), \mathbb{P}_{\theta_{0}})$ with respect to the filtration $\{\sigma(Y_{n}) \colon n \in \natnum\}$ (see Theorem~1, p.442 in \cite{GikhmanSkorokhod}). Given $\varepsilon > 0$, let $B_{\varepsilon}(\theta_{0})$ be a neighborhood of $\theta_{0}$ contained in $\Theta$. We show that with $\mathbb{P}_{\theta_{0}}$ probability one,
\begin{equation} \label{toshow}
    \forall \ M > 0 \ \exists \ N^{\prime} \in \natnum \text{ s.t. } \forall \ \theta \in \Theta \setminus B_{\varepsilon}(\theta_{0}) \ \forall \ n \geq N^{\prime}: \ \operatorname{log}(\varphi_{n}(\theta)) \leq -M. 
\end{equation}
In particular, \eqref{toshow} shows that
\begin{equation} \label{Wald}    
    \mathbb{P}_{\theta_{0}}\Big(\sup\big\{\varphi_{n}(\theta) \, \colon \theta \in \Theta \setminus B_{\varepsilon}(\theta_{0})\big\} \xrightarrow[]{n \to \infty} 0\Big) = 1.
\end{equation}
To show \eqref{toshow}, we let $\theta_{m} \in \partial B_{\varepsilon}(\theta_{0})$ be such that for any $\theta \in \partial B_{\varepsilon}(\theta_{0})$, 
\begin{gather*}
    \operatorname{log}(\varphi_{n}(\theta_{m}))(\omega) \geq \operatorname{log}(\varphi_{n}(\theta))(\omega).
\end{gather*}
This follows from the assumption that for any $\omega \in \Omega$, $\theta \mapsto \varphi_{n}(\theta)(\omega)$ is continuous on $\Theta$. Then, we choose $\omega \in \Omega$ as an element of the set
\begin{equation} \label{oneset}
    \{\omega \in \Omega \colon \forall \ M > 0 \ \exists \ N_{\theta_{m}} \in \natnum \text{ s.t }\forall \ n \geq N_{\theta_{m}}: \ \operatorname{log}(\varphi_{n}(\theta_{m}))(\omega) \leq -M\}.
\end{equation}
By assumption, there exists $N \in \natnum$ such that for any $n \geq N$, $\theta \mapsto \operatorname{log}(\varphi_{n}(\theta))(\omega)$ is concave on $\Theta$. Hence, we put $N^{\prime} \coloneqq \max\{N, N_{\theta_{m}}\}$ and have that $\theta \mapsto \operatorname{log}(\varphi_{n}(\theta))(\omega)$ is concave on $\Theta$ for $n \geq N^{\prime}$. Suppose, by contradiction, that there exists $n \geq N^{\prime}$ and $\theta \in \Theta \setminus B_{\varepsilon}(\theta_{0})$ such that $\operatorname{log}(\varphi_{n}(\theta))(\omega) > -M$. This implies that for $\lambda \in [0,1]$, with $(1-\lambda)\theta_{0} + \lambda \theta \in \partial B_{\varepsilon}(\theta_{0})$,
\begin{align*}
    \operatorname{log}(\varphi_{n}(\theta_{m}))(\omega) &= (1-\lambda)\operatorname{log}(\varphi_{n}(\theta_{m}))(\omega) + \lambda \operatorname{log}(\varphi_{n}(\theta_{m}))(\omega) \\
    &< (1-\lambda)\operatorname{log}(\varphi_{n}(\theta_{0}))(\omega) + \lambda \operatorname{log}(\varphi_{n}(\theta))(\omega), \quad n \geq N^{\prime}, 
\end{align*}
since, for $n \geq N^{\prime}$, $\operatorname{log}(\varphi_{n}(\theta_{0}))(\omega) = 0 > \operatorname{log}(\varphi_{n}(\theta_{m}))(\omega)$ and $\operatorname{log}(\varphi_{n}(\theta_{m}))(\omega) \leq -M$. Then, because of the fact that $\operatorname{log}(\varphi_{n}(\theta))(\omega)$ is concave on $\Theta$ for $n \geq N^{\prime}$, we get,
\begin{gather*}
    \operatorname{log}(\varphi_{n}(\theta_{m}))(\omega) < \operatorname{log}(\varphi_{n}((1-\lambda)\theta_{0} + \lambda \theta))(\omega) \leq \operatorname{log}(\varphi_{n}(\theta_{m}))(\omega), \quad n \geq N^{\prime},
\end{gather*}
which is a contradiction. Thus, since under $\mathbb{P}_{\theta_{0}}$, \eqref{oneset} has measure one, \eqref{toshow} and hence \eqref{Wald} are shown. Following Wald's original consistency proof (Theorems 1 and 2 in \citep{Wald}), we conclude from \eqref{Wald} that $\hat{\theta}_{n} \xrightarrow[]{} \theta_{0}$, as $n \xrightarrow[]{} \infty$, with $\mathbb{P}_{\theta_{0}}$ probability one.
\end{proof}

Given an ML estimator for $\theta_{0}$, we remark that the assumptions on $\varphi_{n}(\theta)(\omega)$ given in Theorem~\ref{consistency} are satisfied if the Hessian matrix of the log-likelihood function $\theta \mapsto \operatorname{log}(p_{n}(\theta))(\omega)$ becomes negative semidefinite for $n$ large enough. As a simple illustration we take $\Theta = [a,b]$, $0< a < b <\infty$, and consider a family $\mathbb{P}_{\sigma^{2}}$, $\sigma^{2} \in [a,b]$, defined upon the exponential family \eqref{exponential}, with known scale parameter $\alpha_{0}$ and unknown variance parameter $\sigma^{2}$. We readily see that the second derivative of the log-likelihood function is negative. Then, if the sequence of coordinates $D$ is such that $D_{+}$ is dense in $\Rplus$, we have seen in Example~\ref{exponentialfam} that the family $\mathbb{P}_{\sigma^{2}}$, $\sigma^{2} \in [a,b]$, is a family of orthogonal Gaussian measures on $\sigma_{D}(\xi)$. In this case, we can apply Theorem~\ref{consistency}, and deduce that variance ML-estimators $(\hat{\sigma}^{2}_{n})$ for the exponential family are strongly consistent. This adds to the results of \cite{Zhang}, in which the sampling domain is assumed to be bounded. In particular, if we remain in the setting of Example~\ref{exampleBM} and consider a $d$-dimensional Brownian motion $B$ starting from zero, we have that $B[\Rplus \cap \mathbb{Q}]_{+}$ is dense in $\Rplus$ with $\mathbb{P}_{B}$ probability one. Thus, if we assume that $\xi$ is sampled along $B[\Rplus \cap \mathbb{Q}]$, a sequence of variance ML-estimators $(\hat{\sigma}^{2}_{n})$ is strongly consistent with $\mathbb{P}_{\scalebox{0.5}{$B$}}$ probability one.